\documentclass{amsart}

\usepackage{amsfonts}
\usepackage{amsmath}
\usepackage{amsthm}
\usepackage{amssymb}
\usepackage[english]{babel}
\usepackage{graphics}
\usepackage{latexsym}
\usepackage{longtable} \usepackage{mathrsfs}
\usepackage{graphicx}
\usepackage{wrapfig} %\usepackage{txfonts}
\usepackage[pdftex,colorlinks=true]{hyperref}

 \newcommand{\larrow}{\longrightarrow}
 \newcommand{\lharpoonup}{-\!\!\!\!\rightharpoonup}

\newcommand{\dist}{\textrm{dist}}

\newcommand{\noi}{\noindent}
\newcommand{\ms}{\medskip}

\newcommand{\La}{\Lambda}
\newcommand{\e}{\varepsilon}
\newcommand{\de}{\delta}

\newcommand{\al}{\alpha}
\newcommand{\ga}{\gamma}
\newcommand{\la}{\lambda}
\newcommand{\p}{\partial}
\newcommand{\Om}{\Omega}

\newcommand{\sub}{\subseteq}
\newcommand{\set}{\setminus}

\newcommand{\R}{\mathbb{R}}
\newcommand{\N}{\mathbb{N}}

\newcommand{\Wga}{[W^{1; \gamma,2}_{\text{aff}}(\R)]^N}

\newcommand{\action}{E (U) = \int_{\R}\left\{ \frac{1}{2}\big| U_x\big|^2 + W \big( U \big) \right\} dx}
\newcommand{\eq}{U_{xx}  = D W ( U )}

\newcommand{\bt}{\begin{theorem}}\newcommand{\et}{\end{theorem}}
\newcommand{\bl}{\begin{lemma}}\newcommand{\el}{\end{lemma}}
\newcommand{\be}{\begin{equation}}\newcommand{\ee}{\end{equation}}
\newcommand{\bc}{\begin{claim}}\newcommand{\ec}{\end{claim}}
\newcommand{\bp}{\begin{proof}}\newcommand{\ep}{\end{proof}}

\newcommand{\BPL}{\medskip \noindent \textbf{Proof of Lemma} }

\newcommand{\BPCO}{\medskip \noindent \textbf{Proof of Corollary} }

\newtheorem{theorem}{Theorem}[section]
\newtheorem{corollary}[theorem]{Corollary}
\newtheorem{lemma}[theorem]{Lemma}

\newtheorem{claim}[theorem]{Claim}

\theoremstyle{definition}

\numberwithin{theorem}{section}

\begin{document}

\title[{On the loss of compactness in the connection problem}]{On the loss of compactness in the vectorial heteroclinic connection problem}

\author{\textsl{Nikos Katzourakis}}
\address{Department of Mathematics and Statistics, University of Reading, Whiteknights Campus, PO Box 220, RG6 6AX, Reading, UK}
%\curraddr{Department of Mathematics and Statistics, Case Western
%Reserve University, Cleveland, Ohio 43403}
\email{n.katzourakis@reading.ac.uk}

%    General info
\subjclass[2010]{Primary 34C37, 46B50, 82B26; Secondary 37K05}

\date{}


\keywords{Heteroclinic connection problem, loss of compactness, phase transitions, Hamiltonian system.}

\begin{abstract} We give an alternative proof of the theorem of Alikakos-Fusco \cite{ AF} concerning existence of heteroclinic solutions $U : \R \larrow \R^N$ to the system
\[ \label{P} \tag{1}
\left\{\begin{array}{l}
  U_{xx} \, =\, D W ( U )\ ,\\
   U(\pm \infty) \, =\, a^{\pm}.\\
\end{array}\right. 
 \]
Here $a^{\pm}$ are local minima of a potential $W \in C^2(\R^N)$
with $W(a^\pm)= 0$. \eqref{P} arises in the theory of phase
transitions. Our method is variational but differs from the original artificial constraint method of \cite{AF} and establishes existence by analysing the loss of compactness in minimising sequences of the action in the appropriate functional space. Our assumptions are slightly different from those considered previously and also imply a priori estimates for the solution. 
\end{abstract}

\maketitle

\section{Introduction.}

In this paper we consider the problem of existence of heteroclinic
solutions to the Hamiltonian ODE system
 \be
\left\{\begin{array}{l} \label{P}
  \eq\ ,\ \ U : \R \larrow \R^N,\\
   U(- \infty) \, =\,  a^{-}\ ,\ \ U(+ \infty) \,=\,  a^{+},\\
\end{array}\right.
  \ee
where $W \in C^2(\R^N)$ is a potential and $a^\pm$ are local minima
of it with $W(a^\pm)=0$. A typical $W$ for $N=2$ is shown in Figures
1,2. Solutions to (\ref{P}) are known as \emph{``heteroclinic
connections''}, being \emph{standing waves}  of the gradient diffusion system
  \be \label{DS} u_t \ = \ u_{xx}
 -   D W(u)\ , \ \ \ \ u\ : \R \times (0,+\infty) \larrow
\R^N.
 \ee
(\ref{P}) arises in the theory of phase transitions. For
details we refer to Alikakos-Bates-Chen \cite{ ABC} and to
Alberti \cite{Al}. From the viewpoint of physics, (\ref{P}) is the Newtonian law of
motion with force $- D (-W)$ induced by the potential $-W$ and
$U$ the trajectory of a test particle which connects two maxima of
$-W$. In the scalar case of $N=1$, existence is textbook material by
phase plane methods. For a variational approach we refer to Alberti
\cite{Al}. Even in this simple case the unboundedness of $\R$ implies that standard compactness and semicontinuity arguments fail when one tries to obtain solutions to $U_{xx}=W'(U)$ variationally as minimisers of the \emph{Action} functional
 \be
\label{A}
 \action.
 \ee
However, for $N=1$ rearrangement methods do apply (Kawohl
\cite{Kaw}). When $N>1$, (\ref{P}) is much more difficult. It has first been
considered by Sternberg in \cite{St}, as a problem arising in the
study of the elliptic system $\Delta U = D W \big(U\big)$. Noting the compactness problems, he
utilises the \emph{Jacobi Principle} to obtain solutions by
studying geodesics in the Riemannian manifold $\big(\R^N \set
\{a^\pm\}, \sqrt{2 W}\langle \_ , \_ \rangle\big)$. 

Following a different approach, Alikakos-Fusco \cite{ AF} subsequently treated \eqref{P}
utilising the \emph{Least Action Principle}. They derived
their solutions as minimisers of \eqref{A}. They introduced an artificial constraint in order to restore
compactness and apply the Direct Method and obtained solutions to
the \eqref{P} by eventually removing the constraint. The same approach has subsequently
been applied by Alikakos jointly with the author \cite{AK} to the respective
travelling wave problem for \eqref{DS}, establishing existence of
solution to the system $\eq - c U_x$ for $c\neq 0$. \eqref{P} has attracted some
attention in connection with the study of system $\Delta U =
D W \big(U\big)$ and related material appears also in
Alama-Bronsard-Gui \cite{ABG}, Bronsard-Gui-Schatzman \cite{BGS}, Alikakos \cite{A,A2} and Alikakos-Fusco \cite{AF3}.

The problem \eqref{P} is nontrivial; except for the failure of the Direct Method for
\eqref{A} due to the loss of compactness, an additional difficulty
when $N>1$ is that the Maximum Principle does not apply. In the papers \cite{AF}, \cite{AK} were introduced substitutes of the Maximum Principle for minimisers. Inspired by these results, the author in \cite{Ka} developed related ideas which apply to general nonconvex functionals. A further difficulty of \eqref{P} is that additional minima of $W$ obstruct existence and suitable assumptions on $W$ must be imposed  (see \cite{AF}).

In the present work, following  \cite{AF}, we obtain
solutions to \eqref{P} as minimisers of \eqref{A}. We bypass their
unilateral constraint method which is of independent interest, but requires a
rather delicate analysis. We establish existence for \eqref{P} by analysing and then restoring by hand the loss of compactness in minimising sequences. Our motivation comes from the theory of
\emph{Concentration Compactness} (see Lions \cite{L1, L2}, and also Bates-Xiaofeng \cite{BX} for a related application of this principle). We note however that Lions' theory merely motivated the ideas utilised herein and we do not know if the well-known condition of ``strict inequality" applies in the present context. Our approach is conceptually different: we introduce a functional space tailored for the
study of \eqref{P} and show that given any minimising sequence
of \eqref{A}, there exist uniformly decaying translates up to which compactness is restored and passage to a minimiser is available (Theorem \ref{Existence - Compactness}). Our
main ingredients are certain energy estimates and measure bounds
which relate to those of \cite{ AF}, \cite{AK}. Herein however we
utilise a different method: we control the behaviour of the minimising
sequence by the sup-level sets $\{W\geq \al\}$ and compactify the
sequence by suitable translations.

Our basic assumption  (A1) is slightly stronger than the
respective of \cite{AF}, but we still allow for a certain degree of
degeneracy. Under this assumption we obtain the a priori quantitative decay estimates
\eqref{*} by means of energy arguments,
without linearising the equation. The rest of the assumptions
 (A2'), (A2'') allow for $W$'s with several minima and
possibly unbounded from below, being similar to those of \cite{AF}. We believe that our proof of the Alikakos-Fusco theorem \cite{AF} provides further insights to the understanding of the problem.

\section{Hypotheses, Setup and the Existence-Compactness Result.}

\noi \textbf{Hypotheses.} We assume $W \in C^2(\R^N)$ with $a^{\pm}$ local minima at zero: $W(a^\pm) =0$. Moreover:

\noi (A1) \emph{There exist $\al_0,\ w_0 >0$ and $\ga \geq 2$ such that for all $\al \in (0,\al_0]$
   the
  sublevel sets $\big\{ W \leq \al \big\}$ contain two $C^2$ stricitly convex components $\big\{ W \leq \al  \big\}^\pm$, each enclosing $a^{\pm}$ respectively such that $\big\{ W = \al \big\}=\p \big\{ W \leq \al \big\}$ and
\[
W(u)\ \geq\ w_0 \big|u\ -\ a^{\pm} \big|^\ga \ , \ \
  \ u \in \big\{W \leq \al_0 \big\}^\pm .
\]}
 
\noindent In addition, \emph{at least one} of the following two properties is satisfied: either

\noi (A2') \emph{we have
\[
\big\{W \leq \al_0 \big\}\ =
    \ \big\{W \leq \al_0 \big\}^+ \bigcup
   \big\{W \leq \al_0 \big\}^- , 
   \]}
or

\noi (A2'') \emph{there exists a convex bounded (localisation) set $\Omega \subseteq \R^N$
    and a $w_{\max} > \al_0$ such that $a^{\pm}$ are global minima of $W\big|_{\Omega}$,
  while
\[
\Om \subseteq  \big\{W \leq w_{\max} \big\},\ \ 
  \p \Om \subseteq \big\{W = w_{\max} \big\}.
\]}
\noi (A1) allows for $C^{\ga-\e}$ flatness at the minima for all $\e>0$ (but not $C^{\infty}$ flatness as in \cite{ AF}, \cite{AK}). The assumption (A2') requires that $\big\{W \leq \al \big\}^\pm$ are the \emph{only} components of the sublevel sets $\big\{W \leq \al \big\}$. We note that there is a crucial local monotonicity assumption hidden inside (A1). this monotonicity is included in the statement that the level sets coincide with the boundaries of the sublevel sets and hence ``flatness" is exluded. 

Under assumption (A2'), we immediately obtain $\liminf_{|u|\rightarrow \infty}W(u) \geq \al_0$. The assumption (A2") allows for $W$'s which may be unbounded from below, assuming nonnegativity of $W$ only within $\Om$.
\[
\underset{\text{Figure 1: A typical } W,\text{ which satisfies assumption (A1) and the
coercivity assumption (A2')}.}{\includegraphics[scale=0.18]{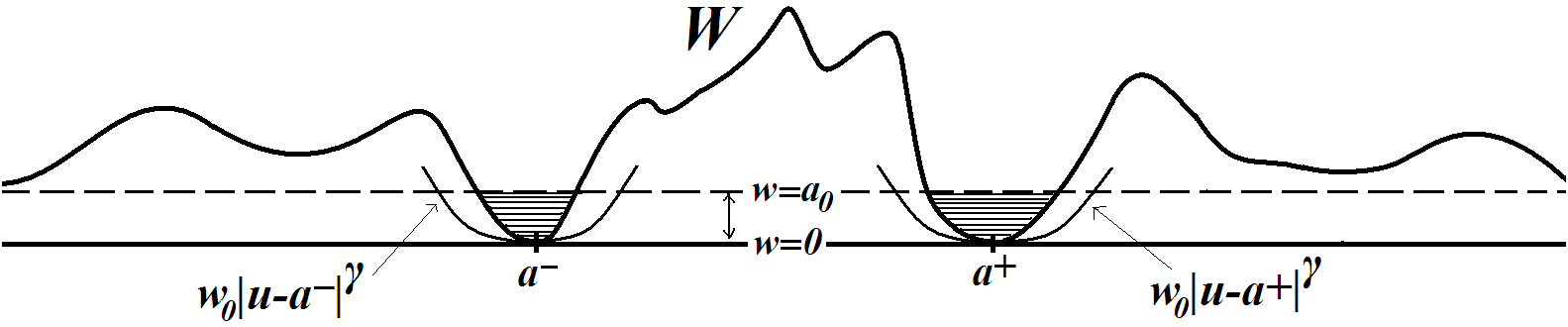}}
\]
\[
\underset{ \text{Figure 2: A typical } W,\text{
the heteroclinic solution } U, \text{ the localisation set } \Om
\text{ of } (A2"), \text{ and the level sets}.}{\includegraphics[scale=0.16]{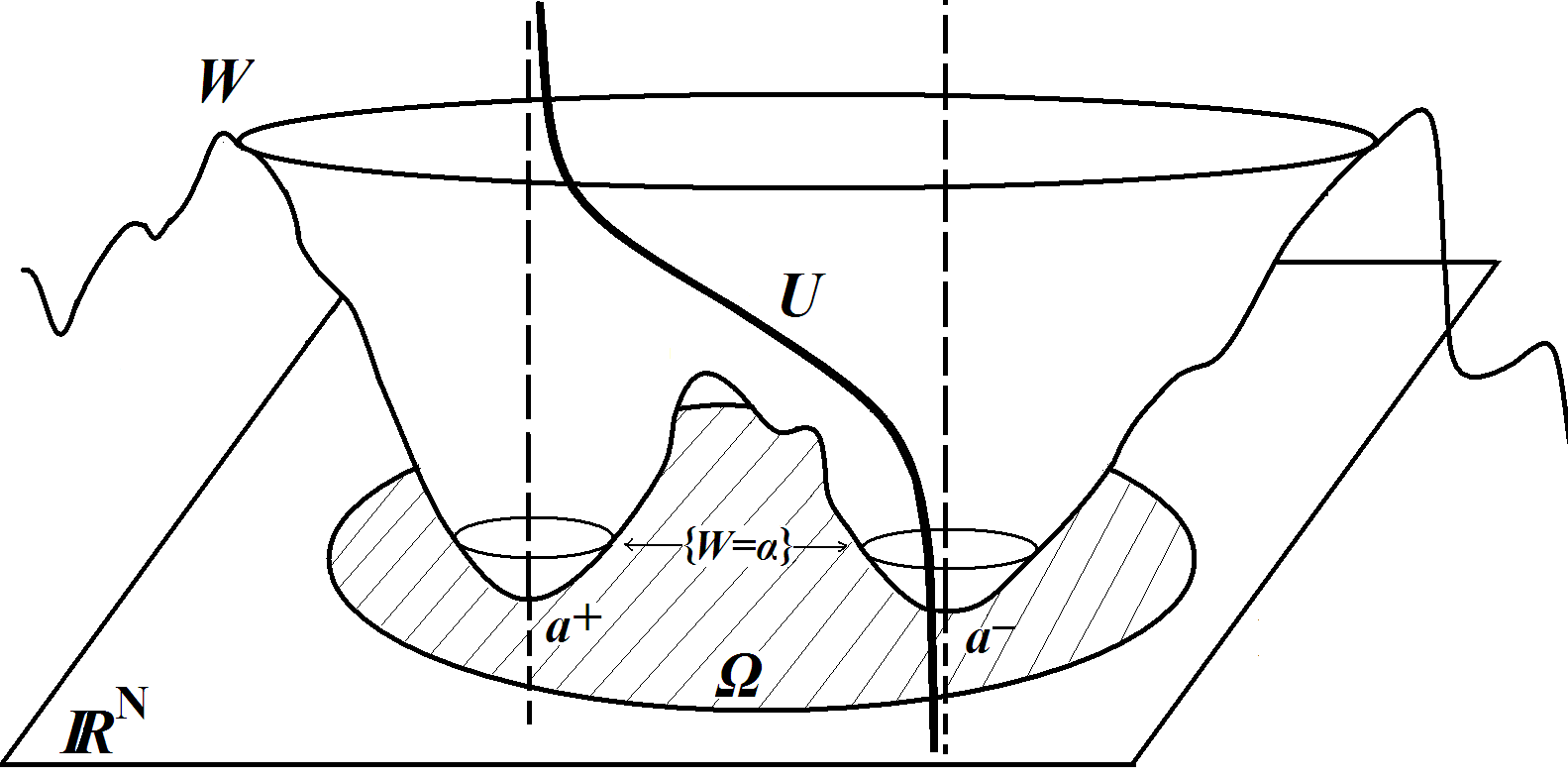}}
\]
\medskip

\noi Under (A2") the existence of a local minimiser $U$ of
\eqref{A} with $E(U)>-\infty$ is a certain issue, but (A1) is
more crucial. We shall refer to (A2') as the
\emph{``coercive''} and to (A2") as the
\emph{``non-coercive''} assumption.

\noi \textbf{Functional setup.} We derive solutions to (\ref{P}) as
minimisers of \eqref{A} in an affine Sobolev space which
incorporates the boundary condition $U(\pm \infty)=a^\pm$ and excludes
the trivial solutions $U=a^\pm$. Let $[W^{1,p}_{\text{loc}}(\R)]^N$
denote the local Sobolev space of vector functions $U\
: \R \larrow \R^N$. For $\e>0$ consider the affine function
 \be \label{eq3}
U_{\text{aff}}^{\e}(x)\ :=\ \left\{
\begin{array}{l}
 a^-\ , \hspace{130pt} x \leq -\e \\
  \left(\dfrac{\e-x}{2\e}\right)a^- \ + \ \left(\dfrac{\e+x}{2\e}\right)a^-\ ,
\ \ \ \  -\e < x <\e\\
 a^+\ , \hspace{130pt}x \geq \e
\end{array}
\right.
 \ee
and set $U_{\text{aff}}^1:=U_{\text{aff}}$. For $p \in
(1,\infty)$, the affine $L^p$-space, $[L^p_{\text{aff}}(\R)]^N :=
[L^p(\R)]^N + U_{\text{aff}}$ is a complete metric space for the $L^p$
distance. The function (\ref{eq3}) will serve also as an a priori upper bound
on the action (\ref{A}) of the minimiser. For $p,\ q \in
(1,\infty)$, we introduce the affine anisotropic Sobolev space
  \be \label{AS} 
[W^{{1;p,q}}_{\text{aff}}(\R)]^N \ := \ \Big\{ U \in [L^p_{\text{aff}}(\R)]^N \ : \ U_x \in [L^q(\R)]^N \Big\}.
 \ee
\eqref{AS} is a complete metric space, isometric to a
reflexive Banach space. The purpose of this work is to establish the following version of the Alikakos-Fusco theorem from \cite{AF}:

\begin{theorem}
 \label{Existence - Compactness} \textbf{(Existence - Compactness)}
Assume that $W$ satisfies (A1) and either (A2') or
(A2"), with $\al_0$, $\ga$, $w_0$, as in (A1),
(A2'), (A2"). There exists a minimising sequence
$(U_i)_1^{{\infty}}$ of the problem 
\[
E(U)= \inf  \Big\{ E(V) : V \in
\Wga \Big\}
\]
for \eqref{A} with $E\big(U_i\big) \geq 0 $. For any
such $(U_i)_1^{{\infty}}$, there exist $(x_i)_1^{\infty} \subseteq
\R$ and translates $\widetilde{U}_i:=U_i(\cdot - x_i)$ which have a
subsequence converging weakly in $\Wga$ to a minimiser $U$ which
solves \eqref{P}:
\[
\left\{\begin{array}{l}
  \eq\ ,\ \ U : \R \larrow \R^N,\\
   U(- \infty) = a^{-}\ ,\ \ U(+ \infty) = a^{+}.\\
\end{array}\right.
\]
In addition, any such minimising solution $U$ satisfies the decay
estimates
 \[ \label{*}  \tag{$*$}
 \left\{\
 \begin{array}{l}
 \big|U(x) - a^\pm \big| \, \leq \,
\left(M{w_0}^{-1}\right)^{\frac{1}{\ga}}
  \ |x |^{-\frac{1}{\ga}},  \ \ \ \ |x| \ \geq \ {M}{\al_0}^{-1},\\
  \hspace{20pt}  \big| U_x(x) \big| \, \leq \, \left(2M\right)^{\frac{1}{2}}\
 |x |^{-\frac{1}{2}},\ \ \ \ \ \ \ \ \ \ |x| \ \geq \ {M}{\al_0}^{-1},
\end{array}\right.
 \]
as well as the bound $E\big(U\big)\leq M$, where
\[
M \,=\, |a^+ - a^- | \underset{[a^-,a^+]}{\max}
\sqrt{2W}.
\]
\end{theorem}

\begin{corollary} \eqref{*} imply that the solution is nontrivial. In
particular, $U\not\equiv a^\pm$.
\end{corollary}

\noi Theorem \ref{Existence - Compactness} asserts that translation
invariance of (\ref{P}) and (\ref{A}) causes the \emph{only}
possible loss of compactness to minimising sequences. The space
$\Wga$ plays a special role to this description. The estimates \eqref{*} are an
essential property, satisfied \emph{uniformly} by the compactified
sequence of the translates and may \emph{not} be satisfied by the
initial $(U_i)_1^{{\infty}}$. In addition they are quantitative, in the sense that the constant depends explicitly on the potential. Moreover, they guarantee that  $U(\pm \infty)=a^\pm$ and
$U_x(\pm \infty)=0$, both fully, not merely up to
subsequences. 

\section{ Proof of the Main Result.}

\noi \textbf{Control on the minimising sequence.} Let $(U_i)_1^{{\infty}}$ be any minimising sequence of (\ref{A}). We will tacitly identify each $U_i$ with its precise representatives. Since
\[
\big|U(x'')-U(x')\big| \ \leq \ (x''-x')^{\frac{1}{2}}\left(
\int_{x'}^{x''}\big|U_x\big|^2dx\right)^{\frac{1}{2}},
\]
we have the inclusion $\Wga \subseteq
[C^{\frac{1}{2}}(\R)]^N$. By (\ref{eq3}), we obtain
\[
E\big(U_{\text{aff}}^{\e}\big)\ = \ \int_{-\e}^\e \left\{
\dfrac{|a^+ - a^-|^2}{8\e^2}\ +\
W\left(\left(\dfrac{\e-x}{2\e}\right)a^- +
\left(\dfrac{\e+x}{2\e}\right)a^- \right)\right\}dx
\]
and hence the explicit bounds
 \be \label{eq5}
\dfrac{\big|a^+ - a^- \big|^2}{4\e}\ \leq \
E\big(U_{\text{aff}}^{\e}\big)\ \leq \ \dfrac{\big|a^+ - a^-
\big|^2}{4\e}\ +\ 2 \e \underset{[a^-,a^+]}{\max}\, W .
 \ee
We immediately get
\[
\inf _{\Wga} E  \ \leq \ \underset{\e>0}{\inf}
\ E\big(U_{\text{aff}}^{\e}\big) \ 
\leq \ \big|a^+ - a^- \big|
\underset{[a^-,a^+]}{\max}\sqrt{2W} \ 
 = \ M \ <\ \infty.
\]
$M$ is necessarily a \emph{strict} upper bound since all
$U_{\text{aff}}^{\e}$ are merely Lipschitz while minimising solutions to
\eqref{P} must be smooth (this latter fact is a consequence of standard regularity considerations of the solutions to the Euler-Lagrange equations). Further, for $i$ large we have
 \be \label{UB}
   \int_{\R}\frac{1}{2} \big| (U_i)_x\big|^2
   dx +  \int_{\R}W\big(U_i\big) dx \ \leq \ M.
 \ee
\noi We now derive $[L^{\infty}(\R)]^N$ bounds. They are
obtained in two different ways, depending on whether (A2') of
(A2") is assumed. In the case of (A2'), it is a
consequence of the next energy estimate. For $\al
\in (0,\al_0]$ and $i=1,\ 2,.. .$ we define the \emph{control
set}
 \be \label{Control set}
\Lambda^\al_i \ := \ \Big\{x \in \R \ : \ W\big(U_i(x)\big) > \al
\Big\}.
 \ee
Let $| \cdot | $ denote the Lebesgue measure on $\R$ and $M$
the constant in estimates \eqref{*}.

\begin{lemma}\textbf{(Energy Estimate I)} \label{Energy Estimate I}
Assume $W$ satisfies (A2'). Then we have
 \begin{align} \label{Estimates1}
M \ &\geq \ \al \ \big| \Lambda^\al_i \big| \ + \ \dfrac{1}{2}
\int_{\R}\big| (U_i)_x\big|^2 dx,\\
 \label{Estimates2}
\big\|{U}_i \big\|_{[L^{\infty}(\R)]^N} \
&\leq \ \big| \Lambda^\al_i \big|^{\frac{1}{2}} \left(\int_{\R}\big|
(U_i)_x\big|^2 dx\right)^{\frac{1}{2}} \ + \ \underset{u \in \{W
\leq \al\}^\pm}{\max}\ |u |,
 \end{align}
 for all $i\in \N$.
\end{lemma}

\BPL \ref{Energy Estimate I}. By  \eqref{UB} and (\ref{Control set}), we have
\begin{align*}
M \ \geq \ E\big(U_i\big)\ & =\; \int_\R W\big(U_i\big)dx
\;+\;\frac{1}{2} \int_{\R}\big| (U_i)_x \big|^2 dx
\\
& \geq \ \int_{\Lambda^\al_i} W\big(U_i\big)dx \;+\;\frac{1}{2}
\int_{\R} \big| (U_i)_x \big|^2 dx
\\
& \geq \ \al\ \big| \Lambda^\al_i \big| \;+\;\frac{1}{2} \int_{\R}
\big| (U_i)_x \big|^2 dx.
 \end{align*}
This proves (\ref{Estimates1}). Let now $(t',t'')$ be
a  subinterval of $\Lambda^\al_i$ such that the endpoints  $U_i(t'),U_i(t'')$ of $U_i\big((t',t'')\big)$ lie on different components of $\{W=\al\big\}$. Hence, we have
\[
\big| U_i(t')-U_i(t'') \big| \ \leq \ \big|t'' -t'\big|^{\frac{1}{2}}
\left(\int_{t'}^{t''} \big| (U_i)_x \big|^2 dx \right)^{\frac{1}{2}} \  \leq \ \big|\Lambda^\al_i \big|^{\frac{1}{2}} \left(\int_{\R} \big|
(U_i)_x \big|^2 dx \right)^{\frac{1}{2}} ,
\]
by using that $U_i(t') \in \big\{W = \al\big\}^\pm$, we
deduce
\[
\big| U_i(t'')-U_i(t') \big|\ \geq \ \big| U_i(t'') \big| \, -\,
\big|U_i(t') \big| \, \geq \ \big| U_i(t'') \big| \, -  \underset{u
\in \{W \leq \al\}^\pm}{\max}\ |u |.
\]
This establishes estimate (\ref{Estimates2}), proving Lemma
\ref{Energy Estimate I}. \qed

\begin{corollary}\textbf{($L^\infty$ bound under (A2'))} \label{Existence of a localised minimising sequence
 in the coercive case} If $W$ satisfies (A1), (A2'), then
 \be
   \label{L infty bounds in the coercive case}
\underset{i \geq 1}{\sup}\ \big\| {U}_i \big\|_{[L^{\infty}(\R)]^N}
\ \leq \ \sqrt{\dfrac{2}{\al_0}}M \ + \, \underset{u \in \{W \leq
\al_0\}^\pm}{\max}\ |u |.
 \ee
\end{corollary}
\noi Now we turn to the case of (A2"). We obtain
existence of a minimising sequence $(U_i)_1^{{\infty}}$ of (\ref{A}) localised inside $\overline{\Om} \subseteq \R^N$ whereon $W\big|_{\Om} \geq 0$. 

\begin{lemma} \label{Localisation}\textbf{($L^\infty$ bound under (A2"))}
\label{Existence of a localised minimising sequence
 in the non coercive case} If $W$ satisfies (A1), (A2"), there is a
minimising sequence $(U_i)_1^{{\infty}}$ for which
$\bigcup_{i=1}^{\infty} U_i(\R) \subseteq  \overline{\Om}$ and
$ W\big(U_i\big) \geq 0$. Moreover,
 \be
  \label{L infty bounds in the non coercive case}
\underset{i \geq 1}{\sup}\ \big\|U_i \big\|_{[L^{\infty}(\R)]^N} \
\leq \  \underset{u \in \p \Om}{\max}\ |u |.
 \ee

\end{lemma}

 \BPL \ref{Existence of a localised minimising sequence
 in the non coercive case}. We show the existence of a  deformation of $W$ to a new
$\overline{W}$ such that $\overline{W}=W$ on $\Om$ and \emph{all}
the minimising sequences of (\ref{A}) relative to
$\overline{W}$ in $\Wga$ can be chosen to be localised inside $\Om$. By (A2"), $W \leq w_{\max}$ inside $\Om$ and $W=w_{\max}$ on $\p \Om$. We define $\overline{W}$ by reflecting with respect to the hyperplane $\{ w=w_{\max} \}$ the portions of the graph of $W$ which lie in the halfspace
$\{w<w_{\max} \}$, to $\{w>w_{\max} \}$. 
\[
\underset{\scriptstyle{\text{Figure 3: The deformed
 coercive potential } \overline{W},\text{ for which } w=w_{\max}
\text{ is a lower bound outside of }
\Om.}}{\includegraphics[scale=0.16]{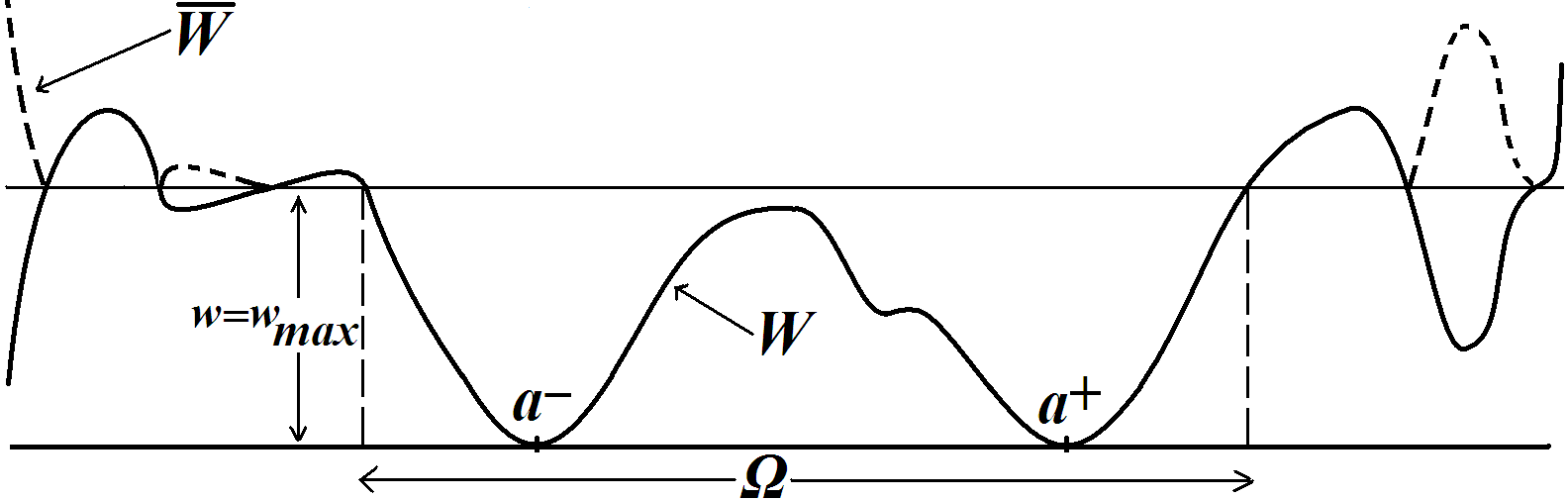}}
\]
By construction, $\overline{W}(u) \geq w_{\max}$, for $u \in \R^N
\set \Om$. Suppose for the shake of contradiction that
$\overline{W}$ has a minimising sequence $(U_i)_1^{{\infty}}$ such
that for some $U_i$ and $a< b$, $U_i\big((a,b)\big)\subseteq \R^N \set \Om$. This is the only case that has to be excluded since by the definition of $\Wga$ the "tails" of each $U_i$ approach asymptotically $a^\pm \in \Om$, at least along a sequence (in general of course there may exists countably many such intervals and we apply this argument to each one of them). By replacing $U_i([a,b])$ by the straight line segment with the same endpoints, i.e. by defining
 \be \label{eq7}
\overline{U}_i(x)\ := \ \left\{\begin{array}{l}
                          U_i(x)\ , \hspace{133pt} x\in \R\set(a,b), \ms\\
                          \left(\dfrac{x-a}{b-a}\right)U_i(b)\ + \
                          \left(\dfrac{b-x}{b-a}\right)U_i(a)\ , \ \ \
                           x \in (a,b),
                        \end{array}
                        \right.
 \ee
we obtain by convexity of $\Om$ that $\overline{U}_i(\R)\subseteq
\overline{\Om}$. By pointwise comparison,
 \be \label{eq8}
\int_a^b \overline{W}(\overline{U}_i(x))dx \ \leq \ \int_a^b
\overline{W}(U_i(x))dx.
 \ee
In addition, $\overline{U}_i\big|_{(a,b)}$ minimises the Dirichlet
integral since it is a straight line, thus
 \be \label{eq9}
\dfrac{\big|\overline{U}_i(b)-\overline{U}_i(a) \big|^2}{b-a}\ = \
\int_a^b \big|(\overline{U}_i)_x \big|^2dx \ < \ \int_a^b
\big|({U}_i)_x \big|^2dx.
 \ee
(\ref{eq8}) and (\ref{eq9}) imply that there exists a minimising sequence of
the Action (\ref{A}) with the potential $\overline{W}$ in the place
of $W$ which lies inside $\overline{\Om}$. Finally,
$W\big|_{\Om}=\overline{W}\big|_{\Om}$ by construction.   \qed

\ms

\noindent In the case that (A2") is assumed, we fix a sequence valued inside $\Om$. Moreover,
\[
M \ \geq \ \underset{i \rightarrow \infty}{\liminf}\ E\big(U_i\big)
\ =: \ \inf \Big\{E(V) : V \in \Wga\Big\} \ \geq \ 0.
\]
As the notation suggests, the right hand side will henceforth
stand for $\underset{i \rightarrow \infty}{\liminf} E\big(U_i\big)$. Now we
employ (A1) to show that $\Lambda^\al_i$ is
connected. For $\al \in (0, \al_0]$, $i=1,2,...$, we set
 \be
  \label{Control times}
   \ \ \ \ \     \la^{\al-}_i \, :=\, \inf \, \Lambda^\al_i \ \ ,\ \ \     \la^{\al +}_i \, :=\, \sup\,
\Lambda^\al_i .
    \ee
We also set
\be \label{da}
d_{\al}\, :=\, \dist\left(\big\{W=\alpha\big\}^-, \, \big\{W=\alpha\big\}^+ \right). 
\ee
We note that $d_\al$ is the distance between the 2 components of the level set $\big\{W=\alpha\big\}$.

\begin{lemma} \textbf{(Control on the $\la^{\al \pm}$ times)}
\label{Control on the times}
 Assume $W$ satisfies (A1) and either (A2') or (A2"). Then, for $\al \in (0,\al_0]$, if $(U_i)_1^\infty$ is the minimising sequence constructed previously, then the respective sets $\Lambda^\al_i$ are intervals and hence
\[ 
\Lambda^\al_i \  =\ \left(\la^{\al -}_i , \, \la^{\al +}_i \right).
\]
\end{lemma}
\BPL \ref{Control on the times}. The claim follows by a direct application of the Replacement Lemma 12 in p.\ 1381 of \cite{AK} by choosing as $\mu$ the Lebesgue measure on $\R$. In order to make the presentation self-contained, we provide also an alternative proof which bypasses this maximum principle type of result of \cite{AK}. We note that the result follows by the replacement lemma of \cite{AF} as well, but this is not entirely direct since herein we use convex level sets and not balls.

We fix a term $U_i$ of the minimising sequence and a respective $\La^\al_i$ and we drop the subscript $i$. Since $\La^\al=\big\{W(U)>\al \big\}$ is open, there exist countably many open intervals such that
\be \label{La}
\La^\al\,=\, \bigcup_{p=0}^\infty \big(x^\al_{2p},x^\al_{2p+1}\big).
\ee
Since $U\in[C^0(\R)]^N$, each image $U\big( \big(x^\al_{2p},x^\al_{2p+1}\big)\big)$ is connected with endpoints on $\big\{W(U)=\al \big\}$ and 
\be \label{eq21}
U(\La^\al)\,=\, \bigcup_{p=0}^\infty U\left( \big(x^\al_{2p},x^\al_{2p+1}\big)\right).
\ee
\begin{claim} \label{cl1} For all $p\in \N$, the image $U\left( \big(x^\al_{2p},x^\al_{2p+1}\big)\right)$ has endpoints on different components $\big\{W(U)=\al \big\}^\pm$ of $\big\{W(U)=\al \big\}$.
\end{claim}
Indeed, supose for the sake of contradiction that for some $p$, both $U\big(x^\al_{2p}\big)$ and $ U\big(x^\al_{2p+1}\big)$ are on $\big\{W(U)=\al \big\}^+$. The deformation of Lemma \ref{Localisation} together with the strictness of assumption (A1) contradicts minimality of $U$. The same holds if the endpoints are on $\big\{W(U)=\al \big\}^-$. The claim follows.

\begin{claim} \label{cl2} 
The set $\La^\al$ consists of finitely many intervals of odd number.
\end{claim}

By Claim \ref{cl1}, for each $p$, $U\left( \big(x^\al_{2p},x^\al_{2p+1}\big)\right)$ has endpoints on different components $\big\{W(U)=\al \big\}$. Hence, in view \eqref{da} we have
\[
d_\al\, \leq\, \big|U\big(x^\al_{2p+1}\big)-U\big(x^\al_{2p}\big)\big|\, \leq \, \int_{x^\al_{2p}}^{x^\al_{2p+1}} |U_x|
\]
and hence for each $q\in\N$, by \eqref{La},
\[
q d_\al\, \leq\, \sum_{p=0}^q \int_{x^\al_{2p}}^{x^\al_{2p+1}}|U_x| \,\leq \, \int_{\La^\al}|U_x|\, \leq \, |\La^\al|^{1/2}\left( \int_{\R}|U_x|^2\right)^{1/2}.
\]
Hence, by Lemma \ref{Energy Estimate I}, we have
\[
q\, \leq\, \frac{1}{d_\al}\left(\frac{M}{\al}\right)^{1/2}M^{1/2}
\]
which implies that there exists a $p_\al \in \N$ no greater than the integer part of $M/\sqrt{\al}d_\al$ such that
\[
\La^\al\,=\, \bigcup_{p=0}^{p_\al} \big(x^\al_{2p},x^\al_{2p+1}\big).
\]
Since
\[
\R\set \La^\al\, =\, \big(-\infty,x^\al_{0}\big ]\bigcup \big[x^\al_{1},x^\al_{2} \big]\bigcup ...\ \ ...\bigcup
\big[x^\al_{2p^{\al}-1}, x^\al_{2p^{\al}} \big]\bigcup \big[x^\al_{2p_{\al}+ 1},+\infty \big)
\]
and $\R\set \La^\al$ equals $\big\{W(U)\leq \al \big\}$, $U$ exits $\big\{W(U)\leq \al \big\}^-$ for the 1st time at $x=x^\al_{0}$ and stays inside $\big\{W(U)\leq \al \big\}^+$ after $x=x^\al_{2p_{\al}+ 1}$  (Figure 4). Since 
\begin{align}
U\big(x^\al_{0}\big) &\in \big\{W=\al \big\}^-, \nonumber\\
 U\big(x^\al_{1}\big), U\big(x^\al_{2}\big)  &\in \big\{W=\al \big\}^+, \nonumber\\
 U\big(x^\al_{3}\big), U\big(x^\al_{4}\big)  &\in \big\{W=\al \big\}^-, \nonumber\\
 & \vdots \nonumber
\end{align}
in view of \eqref{eq21} the number of interval has to odd, for otherwise $U$ stays inside $\big\{W\leq \al \big\}^-$ for infinite time and this contradicts that (at least along a sequence) $U(x)$ converges to $a^+$  as $x\rightarrow \infty$.

\begin{claim} \label{cl3}
All subsets $U\big( \big(x^\al_{1},x^\al_{2}\big) \big)$, $U\big( \big(x^\al_{3},x^\al_{4}\big) \big)$, ... , $U\big( \big(x^\al_{2p_{\al}-1},x^\al_{2p_{\al}}\big) \big)$ of the image $U(\R\set \La^\al)$ lie inside the interior $\big\{W<\al\big\}$ and can not touch the boundary $\big\{W=\al\big\}$ (Figure 4).
\end{claim}

\noi Fix a $q\in \{1,...,p_\al\}$ and assume for the sake of contradiction that there is $[a,b]\sub(x^\al_{2q-1}, x^\al_{2q})$, such that $U([a,b])$ lies on the boundary $\{W=\al\}$. Then, by replacing $U([a,b])$ by the straight line segment with the same endpoints (as in Lemma 3.1), we obtain a contradiction. 

\[
\underset{\scriptstyle{\text{$U(x^{\al_*}_0)$, $U(x^{\al_*}_2)$ can not exist. For brevity we have denoted the points $U(x^\al_p)$ by $x^\al_p$.}}}{\underset{\scriptstyle{\text{Figure 4: Illustration with $p_\al=4$. By minimality the dashed line with endpoints }}}{\includegraphics[scale=0.22]{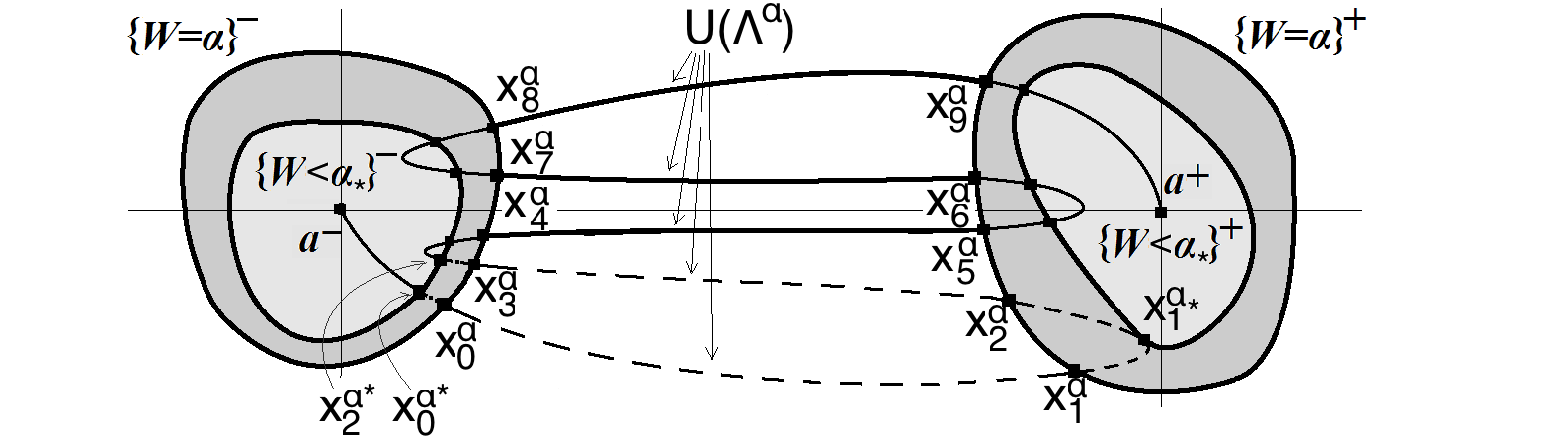}}}
\]

Hence, if $U\big( (x^\al_{2q-1}, x^\al_{2q}) \big)$ touches the boundary $\{W=\al\}$, this happens at isolated points (and otherwise it is inside $\{W<\al\}$). 

Fix such a point and call it $x^*$. By continuity and by assumption (A1), there exist $\de^0,\de^\pm >0$ such that $U\big( (x^*-\de^-, x^*+\de^+) \big)$ lies outside $\{W<\al-\de^0\}$. 
\[
\underset{\scriptstyle{\text{Figure 5.}}}{\includegraphics[scale=0.22]{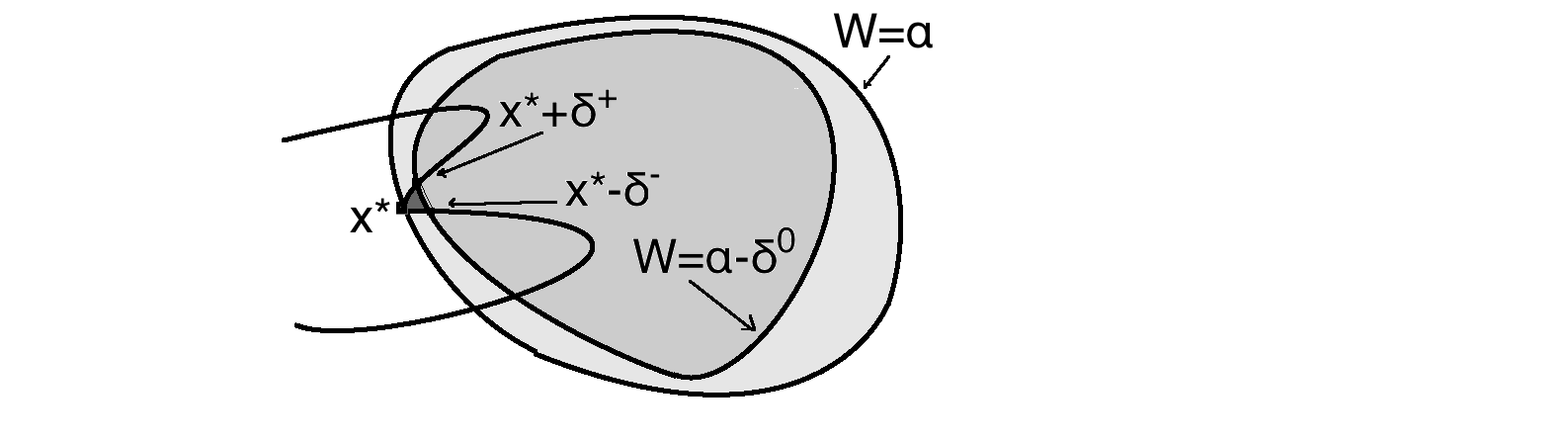}}
\]
By replacing $U\big( (x^*-\de^-, x^*+\de^+) \big)$ by the straight line segment with the same endpoints (as in Lemma 3.1), we obtain a contradiction.  By arguing for all such points $x^*$, we see that 
$U\big( (x^\al_{2q-1}, x^\al_{2q}) \big)$ lies inside $\{W<\al\}$, as desired.

\begin{claim} \label{cl4}
$p_\al =0$, that is $\La^\al$ has only one connected component and hence $x^\al_1=x^\al_{2p_{\al}+1}$.
\end{claim}

We argue by contradiction. Suppose that $p \in \{1,...,p_\al\}$ and consider the set
\be
A\, :=\, \left\{ \beta \in (0,\al)\ \Bigg|\ U\big( \big(x^\al_{2p-1},x^\al_{2p}\big) \big)\bigcap_{p=1}^{p_\al}\big\{ W<\beta \big\}\,  \neq \, \emptyset \right\}.
\ee
Since $U\big( \big(x^\al_{2p-1},x^\al_{2p}\big)$ lies strictly inside the sublevel set, we have that $A\neq \emptyset$. We set
\[
\al_*\, :=\, \inf A.
\]
Since the components $U((x^\al_{2p-1},x^\al_{2p}))$ are finitely many, their distance from the minimum of $W$ is bounded away frow zero and hence $0<\al_*<\al$. By definition of $\al_*$, there exists at least one of the components $U\big( \big(x^\al_{2p-1},x^\al_{2p}\big)\big)$, say for $p=1$, which touches only the boundary of $\big\{ W=\al_* \big\}=\p \big\{ W<\al_* \big\}$ and does not intersect $\big\{ W<\al_* \big\}$. Morover, it can not touch the boundary at more than one points. Hence,
\[
\La^{\al_*}\, =\, \big(x^{\al_*}_{0}, x^{\al_*}_{1}\big)\bigcup \big(x^{\al_*}_{0}, x^{\al_*}_{1}\big)\bigcup\, ...
\]
and consequently $U\big(\big(x^{\al_*}_{0}, x^{\al_*}_{2}\big)\big)$ is contained into $ \big\{ W\geq \al_* \big\}$ and only $U\big(x^{\al_*}_{1}\big)$ is on $ \big\{ W= \al_* \big\}^+$, having both the endpoints $U\big(x^{\al_*}_{0}\big)$, $U\big(x^{\al_*}_{2}\big)$ on $\big\{ W= \al_* \big\}^-$. By arguing as in Lemma 3.1 for $U\big|_{(x^{\al_*}_{0},x^{\al_*}_{2})}$, we obtain a contradiction to the minimality of the action of $U$. Hence, $p_\al=0$.
\ms

By putting Claims  \ref{cl1},  \ref{cl2},   \ref{cl3} and \ref{cl4}, we see that Lemma \ref{Control on the times} has been established.              \qed

\ms

\noi The following sharpens \eqref{Estimates1}, under
the additional information that $\Lambda^\al_i$ is connected.

\bl \textbf{(Energy estimate II)} \label{Energy estimate II}   For all $\al \in (0, \al_0]$ and $i \geq 1$, we have
 \be \label{eq11}
M \ \geq \ E\big(U_i\big) \; \geq \; \dfrac{{d_\al
}^2}{2(\la_i^{\al+}-\la_i^{\al-})} \ + \ \al \
(\la_i^{\al+}-\la_i^{\al-}).
 \ee
 \el

\BPL \ref{Energy estimate II}. Proceeding as
in Lemma \ref{Energy Estimate I}, we recall \eqref{UB} to obtain
\begin{align*}
M \  \geq \ E\big(U_i\big)\ \geq \ \al\ ( \la_i^{\al+} - \la_i^{\al-}) \;+\;\frac{1}{2}
\int_{\la_i^{\al-}}^{\la_i^{\al+}} \big| (U_i)_x \big|^2 dx,
 \end{align*}
 where we have also used Lemma \ref{Control on the
times}. In addition,
\[
d_{\al} \leq \big| U_i(\la_i^{\al-})-U_i(\la_i^{\al+}) \big|
\; \leq \; (\la_i^{\al+}-\la_i^{\al-})^{\frac{1}{2}}
\left(\int_{\la_i^{\al}}^{\la_i^{\al +}} \big| (U_i)_x \big|^2 dx
\right)^{\frac{1}{2}}.
\]
The Lemma follows. \qed

\begin{corollary} \label{Uniform time bounds} \textbf{(Uniform
 bounds on $|\Lambda^\al_i|$)} For $i=1,  2,  . . .$,$\al \in
[0,\al_0]$, we have
 \be \label{eq12}
\dfrac{d_\al^{\ ^2}}{2M}\ \leq \ \big|\Lambda^\al_i \big| \ =\
\la_i^{\al+} \ - \ \la_i^{\al-} \ \leq \ \dfrac{M}{\al}.
 \ee
\end{corollary}

\noi \textbf{Restoration of Compactness.} The bounds (\ref{eq12})
provide information which  allow to control the behaviour of
each $U_i$ by ``tracking''  the $\Lambda^\al_i$'s. In the terminology of \cite{ ABG}, translation
invariance of (\ref{A}) and \eqref{P} allows us to ``fix a centre''
for the $U_i$'s and align the minimising sequence, preventing the terms from escaping to $\pm \infty$. For
$i=1, 2, ...$, we set
 \be \label{eq13}
\hspace{20pt} x_i\ := \ \dfrac{\la_i^{\al_0+} \ + \
\la_i^{\al_0-}}{2}
 \ee
which is the centre of the control set $\Lambda^{\al_0}_i=
\big(\la_i^{\al_0-}, \la_i^{\al_0+}\big)$. We define the
\emph{translates} of the minimising sequence $(U_i)_1^{\infty}$
by:
 \be \label{eq14}
 \hspace{20pt} \widetilde{U}_i \ := \ U_i \big(\cdot - x_i \big)\ , \ \ \ \ i=1,
2, . . .\ .
  \ee
For these translates, their control sets $\widetilde{
\Lambda}^{\al_0}_i = \big(\widetilde{\la}_i^{\al_0 -} ,
\widetilde{\la}_i^{\al_0 +}\big)$ are centred at $x=0$, being
symmetric (Figure 5). The control sets $\widetilde{
\Lambda}^{\al}_i$ of $\widetilde{U}_i$ and $\Lambda^{\al}_i$ of
$U_i$ are related by
 \be \label{eq15}
\big(\widetilde{\la}_i^{\al-},\, \widetilde{\la}_i^{\al+}\big) \ =\
\widetilde{ \Lambda}^{\al}_i \ = \
\left({\la}_i^{\al-}-x_i ,\, {\la}_i^{\al+}-x_i\right).
 \ee
The translates $(\widetilde{U}_i)_1^{\infty}$ defined by
(\ref{eq13}), (\ref{eq14}) will be referred to as the
\emph{compactified sequence} relative to the initial
$(U_i)_1^{\infty}$. The sequence $(\widetilde{U}_i)_1^{\infty}$ will turn out to
be weakly precompact in $\Wga$, converging to a solution of
(\ref{P}).
\[ 
\underset{\text{\!\!\!\!\!\!\!\!\!\!\!\!\!\!\!\!\!\!\!\!\!\!\!\!\!\!\! \!\!\!\!\!\!\!\!\!\!\!\!\!\!\!\! \!\!\!\!\!\!\!\!\!\!\!\!\!\!\!\! \!\!\!\!\!\!\!\!\!
Figure 5: The control sets
of $\widetilde{U}_i$ are symmetric for $\al = \al_0$. For $\al<\al_0$ may not be, but $0 \in \widetilde{\Lambda}^{\al}_i$.} }{\includegraphics[scale=0.21]{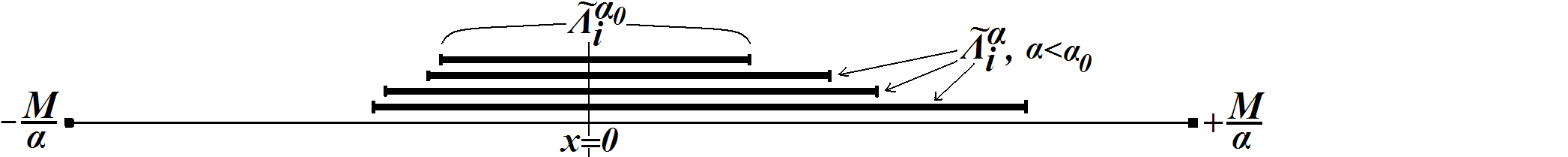} }
\]

\begin{corollary}
  \label{Uniform time bounds for the compactified
sequence} \textbf{(Uniform bounds for the compactified sequence)}
For  $i=1, 2, .. .$ and $\al \in (0,\al_0]$, \eqref{eq12} can be
rewritten in view of \eqref{eq13}, \eqref{eq14}, \eqref{eq15} as
 \be \label{eq17}
\dfrac{d_\al^{\ ^2}}{2M}\ \leq \ \big|\widetilde{\Lambda}_i^{\al}
\big|\ =\ \widetilde{\la}_i^{\al+} \ - \ \widetilde{\la}_i^{\al-} \
\leq \ \dfrac{M}{\al}.
 \ee
In particular, since $0 \in \widetilde{ \Lambda}^{\al}_i$ for 
$\al \in (0,\al_0]$ and $ i=1, 2, . . .$, we have
 \be
\label{Bounds for the translates}
 \medskip \\
 \max \left\{\big| \widetilde{\la}_i^{\al +} \big|
 \ , \  \big| \widetilde{\la}_i^{\al-}\big| \right\} \ \leq \
 \dfrac{M}{\al}.
 \ee
\end{corollary}

\noi \textbf{Bounds and Decay Estimates for the Compactified
Sequence.} The $[L^{2}(\R)]^N$ bound on the derivatives
$(\widetilde{U}_i)_x$ is immediate by the kinetic energy term of
(\ref{A}). The more interesting uniform
$[L^{\gamma}_{\text{aff}}(\R)]^N$ bound is a consequence of our assumption (A1) on the nonconvex
potential term.

\begin{lemma} \label{A priori bounds for the translates}
\textbf{(Estimates for the compactified sequence)}
 Let $(\widetilde{U}_i) _1^{\infty}$ be given by \eqref{eq13} and \eqref{eq14}. If $W$ satisfies
(A1) and either  (A2') or  (A2"), then
$(\widetilde{U}_i) _1^{\infty}$ lies in a ball of $\Wga \cap
[L^{\infty}(\R)]^N$ centred at $U_{\text{aff}}$. Moreover,
\begin{align}
\underset{i \geq 1}{\sup}\ \big\|\widetilde{U}_i -U_{\text{aff}}
\big\|_{[L^{\gamma}(\R)]^N} \ & \leq \
M^{\frac{1}{\gamma}} \left\{\dfrac{1}{w_0}+ \dfrac{2}{\al_0}
\Big\{\underset{i \geq 1}{\sup}\ \big\|\widetilde{U}_i
\big\|_{[L^{\infty}(\R)]^N}\Big\}^\gamma \right\}^{\frac{1}{\gamma}}
\medskip \label{eq18i}
\end{align}
\begin{align}
\underset{i \geq 1}{\sup}\ \big\|\widetilde{U}_i
\big\|_{[L^{\infty}(\R)]^N} \ & \leq \ \left\{\begin{array}{l}
\sqrt{\dfrac{2}{\al_0}}M \ + \ \underset{u
\in \{W \leq \al_0\}^\pm}{\max}\ |u | , \hspace{8pt} \text{ under } (A2')\medskip\\
\medskip
        \underset{u \in \p \Om}{\max} \ |u  |  , \hspace{90pt}
         \text{ under } (A2")
            \end{array}\right.
\medskip \label{eq18ii}
\end{align}
\begin{align}
\underset{i \geq 1}{\sup}\ \big\|\big(\widetilde{U}_i\big)_x
\big\|_{[L^{2}(\R)]^N} \ & \leq \ \sqrt{2M}.
\medskip \label{eq18iii}
\end{align}
\end{lemma}

\BPL \ref{A priori bounds for the translates}.  \eqref{eq18iii}
follows from translation invariance, while \eqref{eq18ii} follows by \eqref{L infty bounds in
the coercive case}, \eqref{L infty bounds in the non coercive case}
and translation invariance. Thus, we only need to prove
\eqref{eq18i}. For,
\[
M\ \geq \ \int_{\R}W\big(U_i\big) dx \ =
 \ \int_{\R}W(\widetilde{U}_i) dx \ \geq \ \int_{-\infty}^{-\frac{M}{\al}}
W(\widetilde{U}_i) dx \
 + \  \int^{+\infty}_{+\frac{M}{\al}}
  W(\widetilde{U}_i) dx\ .
\]
Utilizing (\ref{Bounds for the translates}), we
obtain $W\big(\widetilde{U}_i(x)\big) \leq \al$, for $i=1,2,...$
when $|x | \geq {M}{\al}^{-1}$. Thus, for such $x$ we are in the
domain of validity of (A1). For $\al=\al_0$, we
get
\[
w_0 \left(\int_{-\infty}^{-\frac{M}{\al_0}}\  \Big| \widetilde{U}_i
- a^-\Big|^\gamma dx \ 
 +\  \int^{+\infty}_{+\frac{M}{\al_0}}
 \ \Big| \widetilde{U}_i - a^+\Big|^\gamma dx \right) \ \leq \ M.
\]
By restricting to smaller $\al \leq \al_1 (<\al_0)$, we may assume that
$(-{M}{\al_0}^{-1} , + {M}{\al_0}^{-1}) \supseteq  (-1,1)$. Hence, $U_{\text{aff}}=a^\pm$ for  $|x|\geq {M}{\al_0}^{-1}$. To
conclude, we employ \eqref{eq18ii} to get
\[
\int_{-\frac{M}{\al_0}}^{+\frac{M}{\al_0}} \ \Big| \widetilde{U}_i -
U_{\text{aff}}\Big|^\gamma dx \ \leq \ \dfrac{2M}{\al_0} \left\{
\big\|\widetilde{U}_i \big\|_{ [L^{\infty}(\R)]^N
}\right\}^\gamma.
\]
Putting these estimates together, we see that \eqref{eq18i} has been
established.     \qed

\ms

\begin{lemma} \label{Uniform decay estimate} \textbf{(Uniform decay estimate)}
If $W$ satisfies (A1), the compactified sequence
$(\widetilde{U}_i)_1^{\infty}$ satisfies $ \big|\widetilde{U}_i(x) - a^\pm \big|  \leq 
\left({M}{w_0}^{-1}\right)^{\frac{1}{\ga}}
  \big|x \big|^{-\frac{1}{\ga}}$, for $|x
  |  \geq  {M}{\al_0}^{-1}$.
\end{lemma}

\BPL \ref{Uniform decay estimate}. We have already seen in Lemma
\ref{A priori bounds for the translates} that \eqref{Bounds for the
translates} implies $W\big(\widetilde{U}_i(x)\big) \leq \al$, for
$i=1, 2,... $ when $|x | \geq {M}{\al}^{-1}$. By (A1), we have
\[
w_0\ \big|\widetilde{U}_i(x)-a^{\pm} \big|^\ga  \, \leq\, 
W\big(\widetilde{U}_i(x)\big),
\]
for all such $x \in \R$. Therefore, 
\[
\big|\widetilde{U}_i(x)-a^{\pm}
\big|^\ga  \, \leq\,  \frac{\al}{w_0},
\]
for all $|x | \geq {M}{\al}^{-1}$
and all $\al \leq  \al_0$. We fix an $x \in \R$ for which $|x | \geq
{M}{\al_0}^{-1}$ and choose 
\[
\al(x)\, := \, \frac{|x|}{M}. 
\]
This
is a legitimate choice since $|x|={M} {\al(x)}^{-1} \geq {M}
{\al_0}^{-1}$. We thus obtain that 
\[
\big|\widetilde{U}_i(x)-a^{\pm} \big|^\ga \, \leq\, 
\frac{\al(x)}{w_0}  \leq \, \frac{M}{w_0 |x|}
\]
and by letting $x$ vary,  the estimate follows.
\qed

\ms

\begin{corollary} \label{Decay for solutions}
 \textbf{(A priori decay estimates)}
 Assume $W$ satisfies (A1). Then, if a solution
$U$ to (\ref{P}) exists, it must satisfy estimates
\eqref{*} of Theorem \ref{Existence - Compactness}.
\end{corollary}

\BPCO \ref{Decay for solutions}. We recall from \cite{ AF} the equipartition property $\big|U_x \big|^2  =  2 W\big(U\big)$ satisfied by solutions of (\ref{P}). Equipartition implies
$\big| U_x \big|^2 = 2 W \big( U \big) \leq 2 \al$, for $|x |
\geq {M}{\al}^{-1}$ and $\al \leq \al_0$. The rest follows closely the proof of Lemma
\ref{Uniform decay estimate}. \qed

\ms

\noi \textbf{Passage to a minimising solution.} We conclude by proving
existence of minimisers. By (\ref{eq18i}),
(\ref{eq18ii}) and (\ref{eq18iii}), the sequence of translates
$(\widetilde{U}_i)_1^{\infty}$
converges to some $U$ weakly in $\Wga$ along a subsequence. By denoting the subsequence again by $(\widetilde{U}_i)_1^{\infty}$, we have  that $\widetilde{U}_i - U \lharpoonup 0$  in $[L^{\gamma}(\R)]^N$ and $(\widetilde{U}_i - U)_x \lharpoonup  0$ in $[L^{2}(\R)]^N$, both as $i \rightarrow \infty$. Up to a further subsequence, we have $\widetilde{U}_i \larrow U$ in $[L^{2}_{\text{loc}}(\R)]^N$ and a.e.\ on $\R$ as $i \rightarrow \infty$. By the weak lower semicontinuity of the $L^{2}$ norm and the Fatou
Lemma, we obtain 
\[
E(U) \leq \liminf_{i\rightarrow \infty}E(\widetilde{U}_i).
\]
By (\ref{eq5}), we also get $0\leq E(U) \leq M$. Thus $U$ is a local minimiser of the functional (\ref{A}) in $\Wga$. Hence, $U$ solves (\ref{P}) classically and
satisfies the estimates \eqref{*}. The proof of Theorem \ref{Existence - Compactness} is complete.    \qed

\medskip
\noindent {\bf Acknowledgement.} {We are indebted to N. Alikakos for the careful reading of an earlier version of this manuscript and for providing his valuable suggestions which led to substantial improvements of the content of the paper. We would also like to thank the referee of the paper for their comments and suggestions.}

\end{document}